# Construction of Strong-Uniform Fuzzy Partitions of Arbitrary Dimensions

Zhi Zeng[1,2,a], Ting Wang[1,b]

[a] School of Mechano-Electronic Engineering, Xidian University, 2 South Taibai Road, Xian, Shaan Xi, China

[b] School of Physics and Optoelectronic Engineering, Xidian University, 2 South Taibai Road, Xian, Shaan Xi, China

## Abstract

Strong-uniform fuzzy partition is necessary for the accuracy of fuzzy partition-based histograms. Most previous research focused on constructing one-dimensional strong-uniform fuzzy partitions. While to the best of our knowledge, few have been reported for high-dimensional cases. In order to fill this theoretical vacancy, this paper proves the existence of high-dimensional strong-uniform fuzzy partitions via proposing an analytic formula to construct strong-uniform fuzzy partitions of arbitrary dimensions.

**Keywords: High-dimensional strong-uniform fuzzy partition, membership function, analytic formula**

## 1. Introduction

The histogram is a classical and widely used distribution estimator [1][2]. It is constructed by recording the number of sampling points that fall in each bin. Due to the crisp nature of the histogram, its accuracy is sensitive to the variation of partitions [3][4]. To avoid this problem, people propose using fuzzy partitions to construct histograms, which replace the classical bins with fuzzy subsets [5][6][7]. Among them, the one-dimensional strong-uniform fuzzy partition is the best partitioning [3][8][9]. A lot of theoretical analysis and practical applications have shown that the use of strong-uniform fuzzy partition dramatically increases the accuracy and robustness of the histogram [10][11][12][13][14].

However, the one-dimensional strong-uniform fuzzy partition cannot be used for high-dimensional cases directly. If people want to extend the concept to arbitrarily high dimensions, they should not only explore the existence of high-dimensional strong-uniform fuzzy partitions but also give an effective construction method. Unfortunately, little related research has been reported to the best of our knowledge. This paper proposes an analytic formula to construct a large class of high-dimensional strong-uniform fuzzy partitions to fill this theoretical vacancy.

---

[1] The authors have contributed equally to this paper.
[2] Corresponding author*.
E-mail addresses: zhizeng@mail.xidian.edu.cn (Z. Zeng), tingw2016hg@126.com (T. Wang).

The rest of the paper is organized as follows: Section II introduces the one-dimensional strong-uniform fuzzy partition. Section III proposes the definition of high-dimensional strong-uniform fuzzy partition and an analytic formula to construct a large class of high-dimensional strong-uniform fuzzy partitions. Section IV gives numerical experiments. Section V concludes the paper.

## 2. One-dimensional strong-uniform fuzzy partition

We briefly introduce the one-dimensional strong-uniform fuzzy partition proposed by Kevin Loquin [3]. Take an interval $\Omega=[\mathbf{a}, \mathbf{b}]$ as the universe. Then, the definition of a one-dimensional strong-uniform fuzzy partition is:

**Definition 1:**

Let $m_1 < m_2 < ... < m_p$ be $p$ fixed nodes of the universe, such that $m_1 = a$ and $m_p = b$, and $p \geq 3$. We say that the set of the $p$ fuzzy subsets $A_1, A_2, ..., A_P$, identified with their membership functions $\mu_{A_1}(x), \mu_{A_2}(x), ..., \mu_{A_p}(x)$ defined on the universe, form a strong-uniform fuzzy partition of the universe if they fulfill the following conditions for $k = 1, ..., p$:

1. $\mu_{A_k}(m_k) = 1$ ($m_k$ belongs to what is called the core of $A_k$);

2. If $x \notin [m_{k-1}, m_{k+1}]$, $\mu_{A_k}(x) = 0$ (because of the notation we should add: $m_0 = m_1 = a$ and $m_p = m_{p+1} = b$);

3. $\mu_{A_k}(x)$ is continuous;

4. $\mu_{A_k}(x)$ monotonically increases on $[m_{k-1}, m_k]$ and $\mu_{A_k}(x)$ monotonically decreases on $[m_k, m_{k+1}]$;

5. $\forall x \in \Omega$, $\exists k$, such that $\mu_{A_k}(x) > 0$ (every element of the universe is treated in this partition);

6. $\forall x \in \Omega$, $\sum_{k=1}^{p} \mu_{A_k}(x) = 1$;

7. $\forall k \neq p$, $h_k = m_{k+1} - m_k = h = \text{constant}$, so, $m_k = a + (k-1)h$;

8. $\forall k \neq 2$ and $k \neq p$, $\forall x \in [0, h]$ $\mu_{A_k}(m_k - x) = \mu_{A_k}(m_k + x)$;

9. $\forall k \neq 2$ and $k \neq p$, $\forall x \in [m_k, m_{k+1}]$, $\mu_{A_k}(x) = \mu_{A_{k-1}}(x-h)$ and $\mu_{A_{k+1}}(x) = \mu_{A_k}(x-h)$.

Two typical one-dimensional strong-uniform fuzzy partitions are shown in Fig. 1. In the following, we will explore how to extend the above concepts to the high-dimensional case.

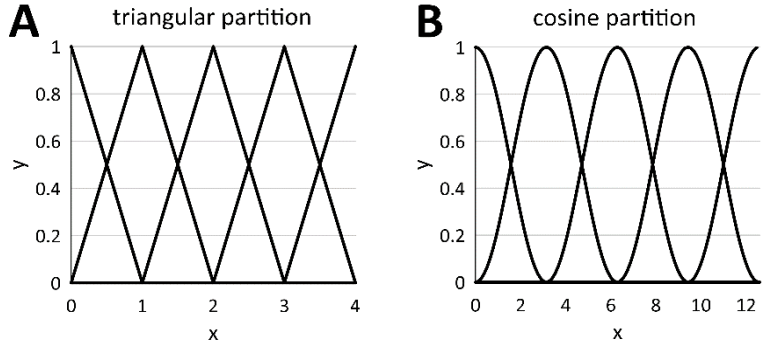

**Figure 1.** One-dimensional uniform strong fuzzy partitions. A. The triangular partition; B. The cosine partition.

## 3. High-dimensional strong-uniform fuzzy partition

This chapter will first propose the definition of high-dimensional strong-uniform fuzzy partition directly. Then, starting from the concept of normalized one-dimensional strong-uniform fuzzy partition, an analytical construction method for high-dimensional strong-uniform fuzzy partition is given.

Generalizing the one-dimensional strong-uniform fuzzy partition to high-dimensional yields the following definition:

**Definition 2:**

Assume that $d$ denotes the number of dimensions. $\forall j \leq d$, let $m_1^j < m_2^j < \cdots < m_{p_j}^j$ be equally spaced nodes in the j$^{\text{th}}$ dimension, that is, $\forall i < p_j$, $m_{i+1}^j - m_i^j = c_j$, and $p_j \geq 3$. The $\prod_{j=1}^{d} p_j$ fuzzy subsets $\{A\}$ (the corresponding set of membership functions is denoted as $\{\mu_A(\mathbf{x})\}$) over the universe $\Omega = [m_1^1, m_{p_1}^1] \times \cdots \times [m_1^d, m_{p_d}^d]$ can form a ***high-dimensional strong-uniform fuzzy partition*** if the membership functions satisfy the following conditions:

1. $\mu(\mathbf{x})$ is non-negative and continuous;

2. Translational symmetry, that is, $\mu_A(\mathbf{x} - \mathbf{x}^A) = \mu_B(\mathbf{x} - \mathbf{x}^B) = \mu(\mathbf{x})$, where $\mathbf{x}^A = \begin{bmatrix} m_{i_1^A}^1 & m_{i_2^A}^2 & \ldots & m_{i_d^A}^d \end{bmatrix}^T$ and $\mathbf{x}^B = \begin{bmatrix} m_{i_1^B}^1 & m_{i_2^B}^2 & \ldots & m_{i_d^B}^d \end{bmatrix}^T$ is the core of fuzzy sets $A$ and $B$. $1 \leq i_j^A, i_j^B \leq p_j$. $\mu(\mathbf{x})$ is the centralized membership function;

3. Mirror symmetry for each dimension. If $\forall i \leq d$, then $|x_i| = |y_i|$, and there must be $\mu(\mathbf{x}) = \mu(\mathbf{y})$;

4. $\mu(\mathbf{0}) = 1$;

5. $\forall \mathbf{x} \in \mathbb{R}^d$, if $\exists 0 \leq k \leq d$ and $|x_j| > c_j$, then $\mu(\mathbf{x}) = 0$;

6. Strong-uniformity, that is, $\forall \mathbf{x} \in \Omega$, $\sum_A \mu_A(\mathbf{x}) = 1$;

7. Radial monotonicity, namely, $\forall \mathbf{x}$ and $\forall \varepsilon > 0$, then $\mu[(1+\varepsilon)\mathbf{x}] \leq \mu(\mathbf{x})$.

It is not difficult to verify that the one-dimensional strong-uniform fuzzy partition is a particular case of **Definition 2**. There are two problems for the high-dimensional and strong-uniform fuzzy partition that needs to be discussed clearly: 1) are there any high-dimensional strong-uniform fuzzy partitions exist; 2) how to construct these partitions. This paper has to define the concept of normalized membership function for one-dimensional strong-uniform fuzzy partitions to answer both questions.

**Definition 3:**

Let $m_1 < m_2 < \cdots < m_p$ be equally spaced nodes on $\mathbb{R}$, i.e., $\forall i < p$ has $m_{i+1} - m_i = c$, and $p \geq 3$. Define the universe as $\Omega = [m_1, m_p]$, $\{A\}$ as the fuzzy partition on it, $\{\mu_A(x)\}$ as the set of corresponding membership functions, and $\mu(x)$ as the corresponding centralized membership function. Then, the ***normalized membership function*** is defined as

$$\eta(x) = \mu(cx) \tag{1}$$

It can be seen from **Definition 3** that, $\eta(x)$ is a function independent of both nodes and intervals. It simply characterizes the shape of the membership function. According to the results of the papers [3][15][16], we know that for any given equally spaced node on $\mathbb{R}$, there always exists infinitely many distinct normalized membership functions for the construction of strong-uniform fuzzy partitions. For the convenience of the subsequent discussion, we denote the set of all normalized membership functions for one-dimensional strong-uniform fuzzy partitions as

$$\mathrm{M} = \{\eta(x)\} \tag{2}$$

The theorem below demonstrates that infinitely many high-dimensional strong-uniform fuzzy partitions can be constructed via $\mathrm{M}$.

**Theorem 2:**

Suppose $d$ denotes the number of dimensions, $\forall j \leq d$, let $m_1^j < m_2^j < \cdots < m_{p_j}^j$ be the equally spaced nodes on the j$^{\text{th}}$ dimension. That is, $\forall i < p_j$ has $m_{i+1}^j - m_i^j = c_j$, and $p_j \geq 3$. Define the universe as $\Omega = [m_1^1, m_{p_1}^1] \times \cdots \times [m_1^d, m_{p_d}^d]$. If one chooses arbitrary $d$ normalized membership functions $\{\eta_j\}_{j=1,2,\cdots d}$ from $\mathrm{M}$, and assume that

$$\mu(\mathbf{x}) = \prod_{j=1}^{d} \eta_j\left(\frac{x_j}{c_j}\right) \tag{3}$$

Then, for any $\mathbf{x}^A = \begin{bmatrix} m_{i_1^A}^1 & m_{i_2^A}^2 & \ldots & m_{i_d^A}^d \end{bmatrix}^T$, and $1 \leq i_j^A, i_j^B \leq p_j$, the fuzzy set $\{A\}$ with $\mathbf{x}^A$ as the fuzzy core and all $\mu_A(\mathbf{x}) = \mu(\mathbf{x} - \mathbf{x}^A)$ act as membership functions constitutes the strong-uniform fuzzy partition on the universe $\Omega$.

**Proof:**

1. The first term of **Definition 2** is obviously satisfied.
2. The second term of **Definition 2** is obviously satisfied
3. For any $\mathbf{x}$ and $\mathbf{y}$, if $|x_i| = |y_i|$ for all $i$, we have

$$\mu(y) = \prod_{j=1}^{d} \eta_j\left(\frac{y_j}{c_j}\right) = \prod_{j=1}^{d} \mu_j(y_j)$$
$$= \prod_{j=1}^{d} \mu_j(x_j) = \prod_{j=1}^{d} \eta_j\left(\frac{x_j}{c_j}\right) \quad (4)$$
$$= \mu(\mathbf{x})$$

Thus, the third term of **Definition 2** is proved.

4. The fourth term of **Definition 2** is satisfied because

$$\mu(\mathbf{0}) = \prod_{j=1}^{d} \eta_j\left(\frac{0}{c_j}\right) = \prod_{j=1}^{d} \mu_j(0) = \prod_{j=1}^{d} 1 = 1 \quad (5)$$

5. The fifth term of **Definition 2** can be verified by substituting the corresponding condition into Eq. (1).

6. As a result of the translational symmetry property of the membership function, it is sufficient to prove that the strong-uniform fuzzy property holds within a bin. Suppose $A$ is the core of a fuzzy set in the fuzzy partition. Any bin $B$ aside of $A$ is a high-dimensional cube with $2^d$ corners. It is helpful to mark these corners as $\{A_{s_1 \cdots s_d}\}$, where $s_j \in \{m_{i_j^A}^j, m_{i_j^A+1}^j\}$. One needs to prove that

$$\sum_{\substack{1 \leq j \leq d \\ s_j \in \{m_{i_j^A}^j, m_{i_j^A+1}^j\}}} \mu_{A_{s_1 \cdots s_d}}(\mathbf{x}) = 1 \quad \forall \mathbf{x} \in B \quad (6)$$

For any positive integer $n < d$, one has

$$\sum_{\substack{1\leq j\leq n+1 \\ s_j\in\{m^j_{i^A_j},m^j_{i^A_j+1}\}}} \mu_{A_{s_1\cdots s_{n+1}}}(\mathbf{x})$$

$$=\sum_{\substack{1\leq j\leq n \\ s_j\in\{m^j_{i^A_j},m^j_{i^A_j+1}\}}} \mu_{A_{s_1\cdots s_n m^{n+1}_{i^A_{n+1}}}}(x_1,\cdots,x_n,x_{n+1})+\mu_{A_{s_1\cdots s_n m^{n+1}_{i^A_{n+1}+1}}}(x_1,\cdots,x_n,x_{n+1})$$

$$=\sum_{\substack{1\leq j\leq n \\ s_j\in\{m^j_{i^A_j},m^j_{i^A_j+1}\}}} \prod_{j=1}^{n}\eta_j\left(\frac{x_j-m^j_{i^A_j}}{c_j}\right)\eta_{n+1}\left(\frac{x_{n+1}-m^{n+1}_{i^A_{n+1}}}{c_{n+1}}\right)+\prod_{j=1}^{n}\eta_j\left(\frac{x_j-m^j_{i^A_j}}{c_j}\right)\eta_{n+1}\left(\frac{x_{n+1}-m^{n+1}_{i^A_{n+1}+1}}{c_{n+1}}\right) \quad (7)$$

$$=\sum_{\substack{1\leq j\leq n \\ s_j\in\{m^j_{i^A_j},m^j_{i^A_j+1}\}}} \prod_{j=1}^{n}\eta_j\left(\frac{x_j-m^j_{i^A_j}}{c_j}\right)\left(\eta_{n+1}\left(\frac{x_{n+1}-m^{n+1}_{i^A_{n+1}}}{c_{n+1}}\right)+\eta_{n+1}\left(\frac{x_{n+1}-m^{n+1}_{i^A_{n+1}+1}}{c_{n+1}}\right)\right)$$

$$=\sum_{\substack{1\leq j\leq n \\ s_j\in\{m^j_{i^A_j},m^j_{i^A_j+1}\}}} \prod_{j=1}^{n}\eta_j\left(\frac{x_j-m^j_{i^A_j}}{c_j}\right)$$

$$=\sum_{\substack{1\leq j\leq n+1 \\ s_j\in\{m^j_{i^A_j},m^j_{i^A_j+1}\}}} \mu_{A_{s_1\cdots s_{n+1}}}(\mathbf{x})$$

In addition, when $n=1$ Eq. (6) holds according to the definition of one-dimensional strong-uniform fuzzy partition. Therefore, the seventh term of **Definition 2** is proved based on mathematical induction.

7. Since

$$\mu[(1+\varepsilon)\mathbf{x}]=\begin{cases}\prod_{j=1}^{d}\eta_j\left(\frac{(1+\varepsilon)x_j}{c_j}\right) & \forall x_j \quad |(1+\varepsilon)x_j|\leq c_j \\ 0 & \text{otherwise}\end{cases}$$

$$\leq\begin{cases}\prod_{j=1}^{d}\eta_j\left(\frac{x_j}{c_j}\right) & \forall x_j \quad |(1+\varepsilon)x_j|\leq c_j \\ 0 & \text{otherwise}\end{cases} \quad (8)$$

$$\leq\begin{cases}\prod_{j=1}^{d}\eta_j\left(\frac{x_j}{c_j}\right) & \forall x_j \quad |x_j|\leq c_j \\ 0 & \text{otherwise}\end{cases}$$

$$=\mu(\mathbf{x})$$

It has radial monotonicity.

∎

So far, this paper has proved the existence of high-dimensional strong-uniform fuzzy partitions and proposed a construction method of a class of fuzzy partitions. It should be noted that the above construction method is not a necessary condition for any strong-uniform fuzzy division. We will demonstrate this conclusion in the next section.

4. **Numerical experiments**

The following section will demonstrate four types of strong-uniform fuzzy partitions: triangular partition, cosine partition, hybrid partition, and a variation of triangular partition.

The two-dimensional triangular partition is shown in Fig. 2. The normalized membership function of the partition is:

$$\eta(\mathbf{x}) = \begin{cases} (1-|x|)(1-|y|) & \|\mathbf{x}\|_\infty \leq 1 \\ 0 & \text{otherwise} \end{cases} \quad (9)$$

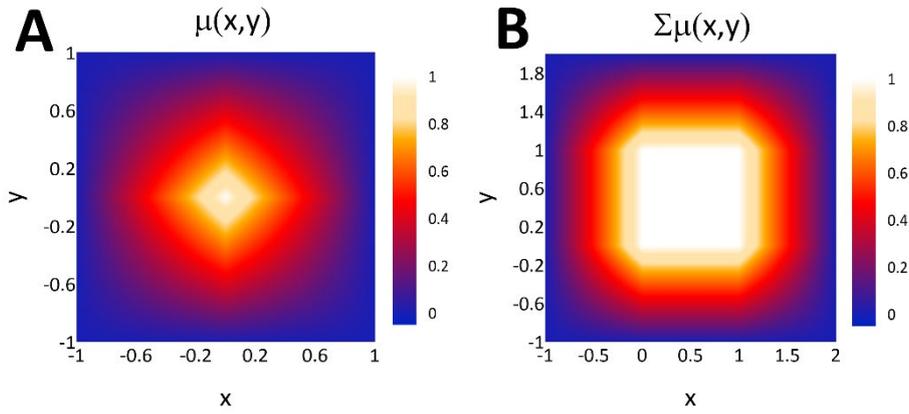

**Figure 2.** Two-dimensional strong-uniform triangular fuzzy partition. A. The centralized membership function; B. The summation of four membership functions around a bin.

Cosine partitions can also be generalized to two-dimensional strong-uniform fuzzy partitions, as shown in Fig. 3. The corresponding normalized membership function is:

$$\eta(\mathbf{x}) = \begin{cases} \dfrac{[\cos(\pi x)+1][\cos(\pi y)+1]}{4} & \|\mathbf{x}\|_\infty \leq 1 \\ 0 & \text{otherwise} \end{cases} \quad (10)$$

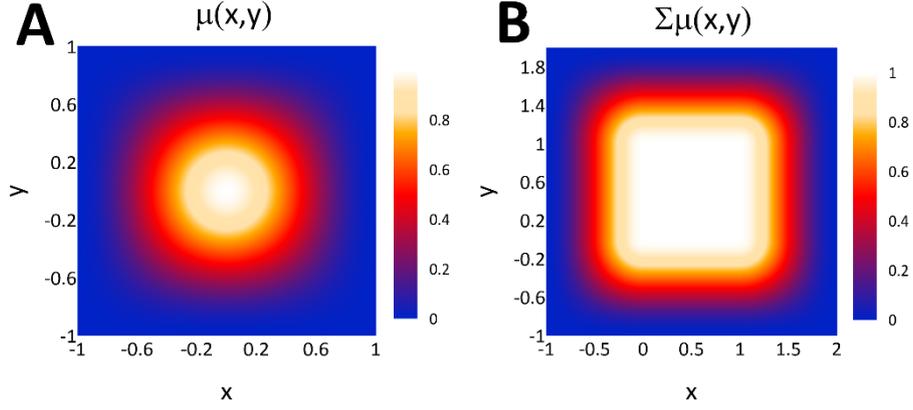

**Figure 3.** Two-dimensional strong-uniform cosine fuzzy partition. A. The centralized membership function; B. The summation of four membership functions around a bin.

It is usually necessary to use the different membership functions for different dimensions. This is because the meanings for different dimensions may be distinct. Fig. 4 shows the case where one dimension adopts a triangular partition, and another adopts a cosine partition. The corresponding normalized membership function is:

$$\eta(\mathbf{x}) = \begin{cases} \dfrac{[\cos(\pi x)+1](1-|y|)}{2} & \|\mathbf{x}\|_\infty \leq 1 \\ 0 & \text{otherwise} \end{cases} \quad (11)$$

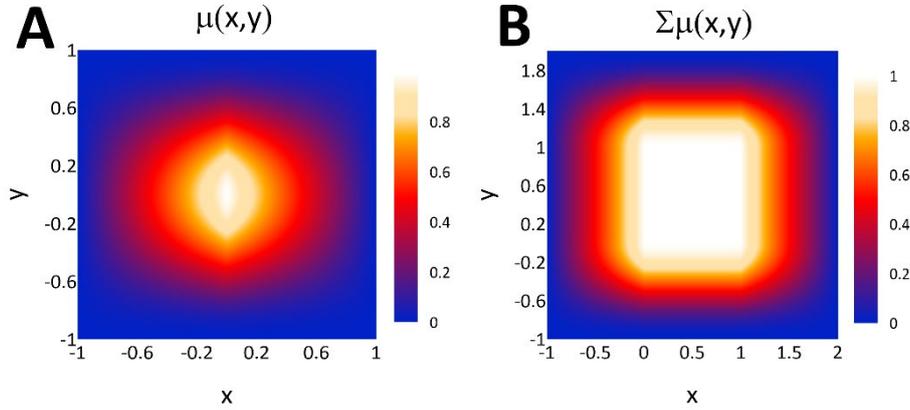

**Figure 4.** Two-dimensional strong-uniform cosine-linear fuzzy partition. A. The centralized membership function; B. The summation of four membership functions around a bin.

It is important to note that the proposed method is a sufficient condition for constructing high-dimensional strong-uniform fuzzy divisions. A variant of the high-dimensional triangular partition can be defined below:

$$\eta(\mathbf{x}) = \frac{1}{2}\left[f_1(\mathbf{x}) + f_2(\mathbf{x})\right] \tag{12}$$

where

$$f_1(\mathbf{x}) = \max\left[0, \min\left(1-|x|, 1-|y|\right)\right] \tag{13}$$

$$f_2(\mathbf{x}) = \max\left[0, \min\left(1-|x-y|, 1-|x+y|\right)\right] \tag{14}$$

It can be seen that the values of Eqs. (9) and (12) are identical when $x=0$ or $y=0$. However, their values are entirely different in other regions. This indicates that the method described in **Theorem 2** is not a necessary condition for constructing high-dimensional strong-uniform fuzzy partitions.

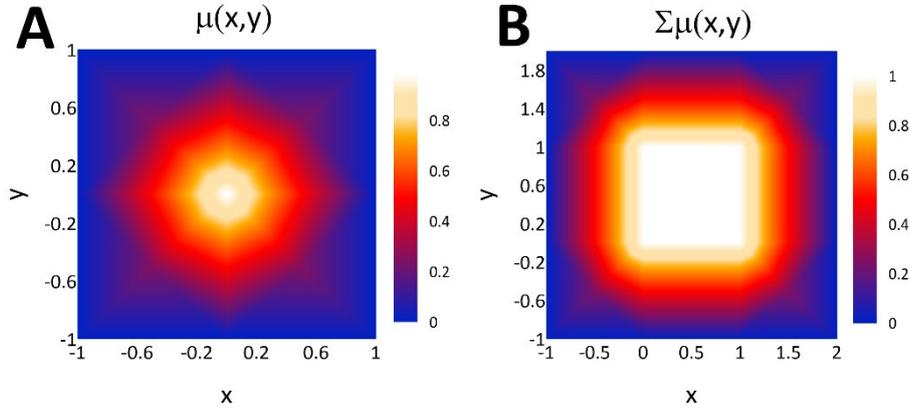

**Figure 5.** A variation of the two-dimensional triangular partition. A. The centralized membership function; B. The summation of four membership functions around a bin.

## 5. Conclusion

This paper first extends the concept of strong-uniform fuzzy partitions to high-dimensional cases. Then, it proposes an analytic formula to construct a large class of high-dimensional strong-uniform fuzzy partitions. The formula indicates that any combination of one-dimensional strong-uniform fuzzy partitions can be used to construct high-dimensional strong-uniform fuzzy partitions. Meanwhile, the formula can also be easily generalized to non-uniform fuzzy partitions. In the subsequent work, we will investigate the necessary conditions for the construction, which would lead to the decomposition theory for generic high-dimensional strong-uniform fuzzy partitions.

*CRediT authorship contribution statement*

**Zhi Zeng:** Conceptualization, Methodology, Formal analysis, Writing; **Ting Wang:** Writing-review & editing, Validation, Experiments.

**Declaration of Competing Interest**

The authors declare that they have no known competing financial interests or personal relationships that could have appeared to influence the work reported in this paper.

**Acknowledgments**

This research is supported by the financial support from the National Natural Science Foundation of China under No. 61805185.